\documentclass[12pt]{article}
\usepackage{amssymb, amsmath}
\begin{document}
\voffset=0.0truein \hoffset=-0.5truein \setlength{\textwidth}{6.0in}
\setlength{\textheight}{8.8in} \setlength{\topmargin}{-0.2in}
\renewcommand{\theequation}{\arabic{section}.\arabic{equation}}
\newtheorem{thm}{Theorem}[section]
\newtheorem{lemma}{Lemma}[section]
\newtheorem{pro}{Proposition}[section]
\newtheorem{cor}{Corollary}[section]
\newcommand{\n}{\nonumber}
\newcommand{\w}{\omega}
\newcommand{\tv}{\tilde{v}}
\newcommand{\tw}{\tilde{\omega}}
\renewcommand{\a}{\alpha}
\renewcommand{\o}{\omega}
\newcommand{\e}{\varepsilon}
\renewcommand{\t}{\theta}
\newcommand{\nao}{\nabla ^\bot \theta}
\newcommand{\vare}{\varepsilon}
\newcommand{\bb}{\begin{equation}}
\newcommand{\ee}{\end{equation}}
\newcommand{\bq}{\begin{eqnarray}}
\newcommand{\eq}{\end{eqnarray}}
\newcommand{\bqn}{\begin{eqnarray*}}
\newcommand{\eqn}{\end{eqnarray*}}
\title{On  the blow-up problem for the axisymmetric 3D Euler equations}
\author{Dongho Chae\\
Department of Mathematics\\
              Sungkyunkwan University\\
               Suwon 440-746, Korea\\
              {\it e-mail : chae@skku.edu}}
 \date{}
\maketitle
\begin{abstract}
In this paper we study the finite time blow-up problem for the axisymmetric 3D incompressible Euler equations with swirl. The evolution equations  for the deformation tensor and the vorticity  are reduced considerably in this case.  Under the assumption of local minima for the pressure on the axis of symmetry with respect to the radial variations we show that the solution blows-up in finite time. If we further assume that the second radial derivative vanishes on the axis, then  system reduces to the form of Constantin-Lax-Majda equations, and can be integrated explicitly.
\end{abstract}

\section{ The axisymmetric  3D Euler equations}
 \setcounter{equation}{0}

 We are  concerned with the following Euler equations for the
homogeneous incompressible fluid flows in a domain $\Omega \subset
\Bbb R^3$,
 \bb\label{e1}
 \frac{D v}{D t} =-\nabla p,
 \ee
 \bb\label{e2}
\textrm{div }\, v =0,
 \ee
 \bb\label{e3}
v(x,0)=v_0(x),
 \ee
  where $D/Dt$ is the material derivative defined by
 $$
 \frac{D}{Dt}=\frac{\partial
 }{\partial t} +(v\cdot\nabla ).
 $$
Here $v=(v_1, v_2, v_3 )$, $v_j =v_j (x, t)$, $j=1, 2,3, $ is the
velocity of the flow, $p=p(x,t)$ is the scalar pressure, and $v_0 $
is the given initial velocity, satisfying div $v_0 =0$.
Since the classical result(\cite{kat}) on the local well-posedness for the
3D Euler equations in the standard Sobolev space $H^m (\Bbb R^3)$, $m>5/2$, the problem of finite time singularity
for such local smooth solution is still an outstanding open problem, although there is a celebrated
result on the blow-up criterion(\cite{bea}) and its refinements(\cite{con1, den, cha0}) taking into account geometric considerations on the vorticity directions.
By an
axisymmetric solution of the Euler equations we mean
velocity field $v(r, x_3, t)$, solving the Euler equations,
and having the representation
$$
v(r, x_3, t)=v^r (r,x_3,t)e_r +v^\theta (r,x_3,t)e_\theta +v^3
(r,x_3,t)e_3
$$
 in the cylindrical coordinate system, where
 $$e_r = (\frac{x_1}{r}, \frac{x_2}{r}, 0), \quad
 e_\theta = (-\frac{x_2}{r}, \frac{x_1}{r}, 0),\quad e_3=(0,0,1),\quad
 r=\sqrt{x_1^2 +x_2^2}.
 $$
In this case also the question of finite time blow-up  of solution
is wide open(see e.g. \cite{cha1,cha2, cha11} for preliminary
studies of the problem; see also \cite{eti} for the related recent result in case of helical symmetry). The vorticity  $\o =$ curl $\,v$ is
computed as
\[ \w = \w^r e_r+ \w^{\t}e_{\t} + \w^3 e_3 ,  \]
where
\[ \w^r = -\partial_{3} v^{\t},
   \quad \w^{\t} = \partial_{3} v^r - \partial_r v^3,
 \quad \w ^3 = \frac{v^\t}{r} +\partial_r v^{\t}.   \]
We denote
\[ \tv = v^r e_r + v^3 e_3 .\]
 The Euler equations for the axisymmetric
solution are
$$
\left\{ \aligned &\partial_t v^r +(\tv\cdot \tilde{\nabla} )v^r
  =-\partial_r p +\frac{(v^\t)^2}{r},\\
  &\partial_t v^3 +(\tv\cdot \tilde{\nabla} )v^3
 =-\partial_{ 3} p,\\
&\partial_t v^\theta +(\tv\cdot \tilde{\nabla}
 )v^\t
 =-\frac{v^r v^\t}{ r} ,\\
 &\textrm{div }\, \tv =0 ,\\
&v(r,x_3,0)=v_0 (r,x_3),
 \endaligned
 \right.
$$
where $\tilde{\nabla} =e_r \partial_ r
  +e_3  \partial_{3}.$
 Note that the above representation of the Euler equations in the cylindrical coordinate system
 is valid off the axis of symmetry, which is chosen to be $x_3-$ axis.
Hence, in order to analyze the equation on the $x_3-$axis, we mainly use the equations in the Cartesian coordinate system.
 Below the functional values of a cylindrically symmetric function $f(x_1,x_2,x_3)=f(r,x_3), r=\sqrt{x_1^2 +x_2^2}$ on the $x_3-$axis, should  be understood as
$\lim_{r\to 0+} f(r,x_3):=f(x_3)$.
\begin{thm}
In the axisymmetric 3D Euler equations with the symmetry of axis chosen to be $x_3-$axis, we write
$$
\bar{\o} :=\partial_1 v^2-\partial_2 v^1(= 2\partial_1 v^2=2\partial_r v^\t), \quad \lambda :=\partial_3 v^3,
$$
which are defined on the $x_3-$axis.
Suppose the initial data satisfies
$$
\mathcal{S}:=\{ x_3 \in \Bbb R \, |\, \bar{\o}_0  (x_3)=0, \lambda_0 (x_3) > 0, \partial_{r}^2p_0(x_3)\geq 0\} \neq \emptyset,
$$
where we denoted
$$ \bar{\o}_0 (x_3)=\bar{\o} (x_3,0),\quad \lambda_0 (x_3)=\lambda(x_3,0), \quad p_0(x_3)=p(x_3,0).$$
 We define $T_1=T_1 (a)$ as
$$T_1  = \inf\{ t>0 \, |\, \partial_{r}^2 p(X_3(a,t), t)<0\},$$
where  $X_3(a,t)$ is the particle trajectory defined by the local classical solution $v(x,t)$.
$$
\frac{\partial X_3(a,t)}{\partial t}= v^3 (X_3(a,t),t), \quad X_3(a,0)=a.
$$
Then, there exists no global classical solution to the axisymmetric 3D Euler equations if
there exists $a\in \mathcal{S}$ such that
\bb\label{1.18}
T_1 (a) \geq  \frac{1}{\lambda_0 (a)}.
\ee
\end{thm}
{\bf Remark 1.1 } After the above theorem is proved P. Constantin informed me of the preprint(\cite{con0}), where it is shown that the positivity of all the matrix components of the hessian of the pressure leads to a singularity in the general case. In the above, however, the positivity is assumed essentially only for one component on the $x_3-$axis in the case of axisymmetry.\\
\ \\
\noindent{\bf Remark 1.2 } The assumption on the positivity of the second radial derivative of the pressure on the $x_3-$axis is physically natural in view of the following heuristic argument. We consider a axisymmetric compressible  ideal fluid with swirl. Due to centrifugal force the density of fluid becomes local minimum on the $x_3-$axis, which implies local minimum of pressure on the $x_3-$axis, hence $\partial _r ^2 p\geq 0$ on the axis. Now we take zero Mach number limit for the pressure to obtain the pressure of the original axisymmetric incompressible fluid(see \cite{kla} for rigorous result for this singular limit problem). In this limiting procedure it is  plausible to expect preservation of the local minimum property ofthe pressure on the $x_3-$axis.\\
\ \\
\begin{thm}
In the axisymmetric 3D Euler equations with the symmetry of axis chosen to be $x_3-$axis,
let us assume that  there exists $T>0$ such that
\bb\label{1.18}
\partial_r ^2 p(x_3,t)=0 \quad \forall (x_3, t)\in \Bbb R \times [0, T].
\ee
Then, the pair $(\bar{\o}, \lambda)$, which is defined in Theorem 1.1, can be explicitly given by
\bq\label{1.19}
\bar{\o }(X_3(a,t),t)&=& \frac{4\bar{\o}_0 (a)}{(2-\lambda _0 (a) t)^2 +\bar{\o}_0 (a) ^2 t},\\
\label{1.19a}
\lambda (X_3(a,t),t)&=&\frac{4\lambda_0 (a)-2[\lambda_0 (a) ^2 +\bar{\o}_0 (a)^2 ]t}{
(2-\lambda _0 (a) t)^2 +\bar{\o}_0 (a) ^2 t}
\eq
along the particle trajectory $\{ X_3 (a,t)\}$ for all $(a, t)\in \Bbb R\times [0, T]$.
Let us assume $\mathcal{S}_0=\{ x_3 \in \Bbb \Bbb \, |\, \lambda_0 (x_3 )>0,\, \bar{\o}_0 (x_3 )=0\}\neq \emptyset$.  Then,  the form of solution (\ref{1.19a}) implies that there exists no global classical solution to the 3D axisymmetric Euler equations  if
$$T\geq \inf_{a\in \mathcal{S}_0} \frac{2}{\lambda _0 (a) }.
$$
\end{thm}

\section{Proof of the main theorems}
 \setcounter{equation}{0}

We begin with the following elementary lemma.
\begin{lemma}
Let $v=(v^1,v^2,v^3)=v(x_1,x_2,x_3)$ be an axially symmetric $C^1-$vector field a on $\Bbb R^3$ with the axis of symmetry chosen as the $x_3-$axis, satisfying  div $v=0$, and let $p=p(x_1,x_2,x_3)$ be an axially symmetric $C^2-$scalar function on $\Bbb R^3$. Then,
on the axis of symmetry we have
\bq\label{1.10a}
&&v^1=v^2=\partial_3 v^1=\partial_3 v^2=\partial_1 v^3=\partial_2 v^3=0,  \\
\label{1.10ab}
&&v^r=v^\t=\partial_3 v^r=\partial_3 v^\t=\partial_r v^3=0,
\eq
\bq\label{1.11}
&&\partial_1 v^1 = \partial_2 v^2=-\frac{\partial_3 v^3}{2}=\partial_r v^r=\lim_{r\to 0}\frac{v^r}{r},\\
\label{1.11aa}
&&\partial_1 v^2=-\partial_2 v^1=\partial_r v^\t =\lim_{r\to 0}\frac{v^\t}{r},
\eq
\bb\label{1.11a}
\partial_1 p=\partial_2 p=\partial_1\partial_2 p=\partial_1\partial_3 p=\partial_2\partial_3 p=\partial_r p=\partial_r\partial_3 p=0,
\ee
\bb\label{1.11ab}
\quad \partial _1^2 p=\partial_2^2 p=\partial_r^2 p=\lim_{r\to 0}\frac{\partial_r p}{r}.
\ee
\end{lemma}
 {\bf Proof }
Here we use notations,
$$(x_1^\prime, x_2^\prime):=(-x_2,x_1)\quad \mbox{and}\quad (\bar{x}_1, \bar{x}_2):=(\frac{r}{\sqrt{2}}, \frac{r}{\sqrt{2}}), \quad r=\sqrt{x_1^2+x_2^2}.
$$
Let us observe first
$$
v^1 (x_1,x_2, x_3)= \frac{x_1}{r} v^r -\frac{x_2}{r} v^\t , \quad v^2(x_1,x_2, x_3)= \frac{x_2}{r} v^r+\frac{x_1}{r} v^\t, \quad r>0,
$$
and thus
\bq\label{le1}
v^1 (x_1^\prime,x_2^\prime, x_3)&=& -\frac{x_2}{r} v^r -\frac{x_1}{r} v^\t =-v^2(x_1,x_2, x_3),\\
\label{le2}
v^2(x_1^\prime,x_2^\prime, x_3)&= &\frac{x_1}{r} v^r-\frac{x_2}{r} v^\t =v^1 (x_1,x_2, x_3).
\eq
Passing $r\to 0$ in (\ref{le1})-(\ref{le2}), we find that
$v_1=v_2=0$ on the $x_3-$axis.   On the other hand,
\bq
\label{le3}
v^1 (\bar{x}_1,\bar{x}_2, x_3)+v^2 (\bar{x}_1,\bar{x}_2, x_3)&=& \sqrt{2} v^r(x_1,x_2, x_3) ,\\
\label{le4}
 v^1 (\bar{x}_1,\bar{x}_2, x_3)-v^2 (\bar{x}_1,\bar{x}_2, x_3)&=& -\sqrt{2} v^\t(x_1,x_2, x_3) ,
\eq
and passing $r\to 0$ in (\ref{le3})-(\ref{le4}), we also find that $v^\t=v^r=0$ on the $x_3-$axis.
Replacing $v_1, v_2, v^r,v^\t$ by $\partial_3 v^1, \partial_3 v^2, \partial_3 v^r,\partial_3 v^\t$ respectively in the above argument we also deduce that
$\partial_3 v^1=\partial_3 v^2=\partial_3 v^r=\partial_3 v^\t=0$. Next we note that
$$
\partial_1 v^3(x_1,x_2, x_3)= \frac{x_1}{r}\partial_r v^3, \quad
\partial_2 v^3(x_1,x_2, x_3)= \frac{x_2}{r}\partial_r v^3,\quad r>0,
$$
and therefore
\bq
\label{le5}
\partial_1 v^3(x_1^\prime,x_2^\prime, x_3)&=&-\frac{x_2}{r}\partial_r  v^3=-\partial_2 v^3(x_1,x_2, x_3), \\
\label{le6}
\partial_2 v^3(x_1^\prime,x_2^\prime, x_3)&=& \frac{x_1}{r}\partial_r v^3=\partial_1 v^3(x_1,x_2, x_3)
\eq
for all $r>0$.
Similarly to the above, passing $r\to 0$ in (\ref{le5})-(\ref{le6}), we deduce
$\partial_1 v^3=\partial_2 v^3=0$ on the $x_3-$axis.  Since
$$
\partial_1 v^3(\bar{x}_1,\bar{x}_2, x_3)-\partial_2 v^3(\bar{x}_1,\bar{x}_2, x_3)=\sqrt{2} \partial_r v^3(x_1,x_2, x_3),
$$
we are lead to $\partial_r v^3=0$ on the $x_3-$axis by passing $r\to 0$.
In order to verify (\ref{1.11}) we compute
\bqn
\partial_1 v^1 (x_1,x_2, x_3)&=&\frac{v^r}{r}-\frac{x_1^2}{r^3} v^r +\frac{x_1^2}{r^2}\partial_r v^r +\frac{x_1 x_2}{r^3} v^\t -\frac{x_1 x_2}{r^2}\partial_r v^\t,\\
\partial_2 v^2(x_1,x_2, x_3)&=&\frac{v^r}{r}-\frac{x_2^2}{r^3} v^r +\frac{x_2^2}{r^2}\partial_r v^r -\frac{x_1 x_2}{r^3} v^\t +\frac{x_1 x_2}{r^2}\partial_r v^\t
\eqn
for $r>0$.
Hence,
\bb
\label{le7}
\partial_1 v^1 (x_1^\prime,x_2^\prime, x_3)=\partial_2 v^2(x_1,x_2, x_3)\quad \forall r>0.
\ee
Passing $r\to 0$ in (\ref{le7}), we have
$\partial_1 v^1=\partial_2 v^2$ on the $x_3-$axis. The condition div $v=0$ implies $\partial_1 v^1=\partial_2 v^2=-\frac12 \partial_3 v^3$.
We note
\bb\label{1.10aa}
\partial_1 v^1 (\bar{x}_1,\bar{x}_2, x_3)-\partial_2 v^2 (\bar{x}_1,\bar{x}_2, x_3)=\frac{v^\t(x_1,x_2, x_3)}{r} -\partial_r v^\t(x_1,x_2, x_3),
\ee
and
\bb
\label{1.10ab}
\partial_1 v^1 (\bar{x}_1,\bar{x}_2, x_3)+\partial_2 v^2 (\bar{x}_1,\bar{x}_2, x_3)=\frac{v^r(x_1,x_2, x_3)}{r} +\partial_r v^r(x_1,x_2, x_3).
\ee
From  (\ref{1.10aa}) we have
\bb\label{1.10aba}
\lim_{r\to 0}\frac{v^\t}{r} =\lim_{r\to 0}\partial_r v^\t.
\ee
Let us compute
\bqn
\partial_2 v^1 (x_1,x_2, x_3)&=&-\frac{x_1 x_2}{r^3}v^r +\frac{x_1 x_2}{r^2}\partial_r v^r -\frac{v^\t}{r} +\frac{ x_2^2}{r^3} v^\t -\frac{ x_2^2}{r^2}\partial_r v^\t,\\
\partial_1 v^2(x_1,x_2, x_3)&=&-\frac{x_1 x_2}{r^3}v^r +\frac{x_1 x_2}{r^2}\partial_r v^r +\frac{v^\t}{r} -\frac{x_1 ^2}{r^3} v^\t +\frac{ x_1^2}{r^2}\partial_r v^\t,
\eqn
and, hence
\bb\label{le8}
\partial_1 v^2 (x_1^\prime,x_2^\prime, x_3)=-\partial_2 v^1(x_1,x_2, x_3)
\ee
for all $r>0$.
Passing $r\to 0$ in (\ref{le8}), we obtain
$\partial_1 v^2=-\partial_2 v^1$ on the $x_3-$axis. Let us compute
\bb\label{1.10ac}
\partial_2 v^1 (\bar{x}_1,\bar{x}_2, x_3) +\partial_1 v^2(\bar{x}_1,\bar{x}_2, x_3)=-\frac{v^r(x_1,x_2, x_3)}{r} +\partial_r v^r(x_1,x_2, x_3),
\ee
and
\bb
\label{1.10ad}
\partial_2 v^1 (\bar{x}_1,\bar{x}_2, x_3) -\partial_1 v^2(\bar{x}_1,\bar{x}_2, x_3)=-\frac{v^\t(x_1,x_2, x_3)}{r} -\partial_r v^\t(x_1,x_2, x_3).
\ee
The equation (\ref{1.10ac}) provides us with
\bb\label{1.10ae}
\lim_{r\to 0}\frac{v^r}{r}=\lim_{r\to 0}\partial_r v^r,
\ee
 while (\ref{1.10ad}), combined with (\ref{1.10aba}),
shows $\partial_1 v^2= \partial_r v^\t$ respectively on the $x_3-$axis.
Using the fact (\ref{1.10ae}), passing $r\to 0$ in (\ref{1.10ab}), we deuce
$\partial_1 v_1 =\partial_r v_r$ on the $x_3-$axis.
As for (\ref{1.11a}) the proof of $\partial_1 p=\partial_2 p=\partial_1\partial_3 p=\partial_2\partial_3 p=\partial_r p=\partial_r\partial_3 p=0$ is exactly same as the above and we omit it. We note
$$
\partial_1\partial_2 p(x_1,x_2, x_3)=-\frac{x_1 x_2}{r^3} \partial_r p + \frac{x_1 x_2}{r^2} \partial_r ^2 p,
$$
and find that $\partial_1\partial_2 p(x^\prime_1,x^\prime_2, x_3)=-\partial_1\partial_2 p(x_1,x_2, x_3)$.
Hence, passing $r\to 0$, we have $\partial_1\partial_2 p=0$ on the $x_3-$axis.
We also compute
$$
\partial_1^2  p(x_1,x_2, x_3)=\frac{1}{r}\partial_r p-\frac{x_1^2}{r^3} \partial_r p +\frac{x_1^2}{r^2}\partial_r ^2 p=\partial_2^2  p(x^\prime_1,x^\prime_2, x_3),
$$
and deduce that
\bb\label{1.11aa}
\partial_1^2  p =\partial_2^2  p \quad \mbox{on the $x_3-$axis by passing $r\to 0$}
\ee
Note that
\bb\label{1.11b}
\partial_1^2  p +\partial_2^2  p =\frac{1}{r}\partial_r p +\partial_r ^2 p,
\ee
and
\bb\label{1.11c}
\partial_1\partial_2 p(\bar{x}_1,\bar{x}_2, x_3)=-\frac{1}{2r} \partial_r p + \frac12 \partial_r ^2 p\to 0 \quad\mbox{as $r\to 0$}.
\ee
From (\ref{1.11aa}), (\ref{1.11b}) and (\ref{1.11c}) we have
$\partial _1^2 p=\partial_2^2 p=\partial_r^2 p=\lim_{r\to 0}\frac{1}{r}\partial_r p$.
$\square$\\
\ \\
Next we recall the matrix representation of the Euler equations(see e.g. \cite{maj}).
Given velocity $v(x,t)$, and pressure $p(x,t)$, we introduce the
$3\times 3$ matrices,
$$
V_{ij}=\frac{\partial v_j}{\partial x_i},\quad
S_{ij}=\frac{V_{ij}+V_{ji}}{2},\quad A_{ij}=\frac{V_{ij}-V_{ji}}{2},
\quad P_{ij}=\frac{\partial ^2 p}{\partial x_i \partial x_j},
$$
with  $i,j=1,2,3$. Then,  we have the decomposition
$V=(V_{ij})=S+A$, where $S=(S_{ij})$ represents the deformation
tensor of the fluid, and $A=(A_{ij})$ is related to the vorticity
$\o$
 by the formula,
 \bb \label{1.4}
 A_{ij}= \frac12 \sum_{k=1}^3 \e
_{ijk} \o_k,\qquad \o_i = \sum_{j,k=1}^3\e_{ijk}A_{jk},
 \ee
 where $\e_{ijk}$ is the skewsymmetric tensor with the normalization
 $\e_{123}=1$. Note that $P=(P_{ij})$ is the hessian of the
 pressure.
  Let  $\{ \lambda _1 , \lambda _2 , \lambda _3\}$
be the set of eigenvalues  of $S$. Computing partial derivatives
$\partial/
\partial x_k$ of (\ref{e1}) yields
 \bb\label{1.5}
 \frac{D V}{Dt}=-V^2 -P.
\ee
 Taking symmetric part of (\ref{1.5}), we have
 \bb\label{1.6}
 \frac{D S}{Dt}=-S^2-A^2-P,
 \ee
 from which, using the formula (\ref{1.4}), we derive
 \bb\label{1.7}
 \frac{D S_{ij}}{Dt} = -\sum_{k=1}^3 S_{ik}S_{kj}
 +\frac14 (|\o|^2 \delta_{ij} -\o_i\o_j )- P_{ij},
 \ee
 where $\delta_{ij}$ is the Kronecker delta defined by
 $\delta_{ij}=1$ if $i=j$, and $\delta_{ij}=0$ otherwise.
 The antisymmetric part of (\ref{1.5}) is
 \bb\label{1.8}
 \frac{DA}{Dt}=-SA -AS,
 \ee
 which, using the formula (\ref{1.4}) again, we obtain easily
 \bb\label{1.9}
 \frac{D \o}{Dt} = S \o,
 \ee
 which is the well-known vorticity evolution equation that could be
 derived
  also by taking curl of (\ref{e1}).
  Taking trace of (\ref{1.7}), we have the identity
  \bb\label{1.10}
  -(\lambda_1 ^2 +\lambda_2 ^2 +\lambda_3 ^2) + \frac12 |\o|^2=\Delta p.
  \ee

\noindent{\bf Proof of Theorem 1.1 }
 Thanks to lemma 2.1 we have the following reduced representation for the deformation tensor, the vorticity, and the hessian of the pressure  on the $x_3-$axis.
  \bb\label{1.12}
   S=\mathrm{diag}(-\frac{\lambda}{2}, -\frac{\lambda}{2} , \lambda ),\quad
   \o=(0,0, \bar{\o }),\quad
 P=\mathrm{diag}(\partial_{r} ^2 p,\partial_{r} ^2 p, \partial_{3} ^2 p),
\ee
where
$$ \lambda=\partial_3 v_3\left(=-\frac12 \partial_1 v^1=-\frac12 \partial_2 v^2=-\frac12 \partial_r v^r\right), \quad \bar{\o}=\partial_1 v^2-\partial_2 v^1\left(=2\partial_r v^\t\right)
$$
on the $x_3-$axis.
The $(11)$ and $(22)$ components of the matrix equation (\ref{1.7}) reduce to
 \bb\label{1.14}
\frac{\bar{D }\lambda  }{Dt}=\frac{\lambda^2}{2}-\frac{\bar{\o} ^2}{2}+2\partial_{r}^2 p,
\ee
where we set
$$ \frac{\bar{D}}{Dt}=\partial_t + v_3\partial_3 ,
$$
while the $(33)$ component becomes
\bb\label{1.15}
\frac{\bar{D }\lambda  }{Dt}= -\lambda^2 -\partial_{3}^2 p.
  \ee
We note that (\ref{1.10}) reduces to
  \bb\label{1.16}
  \Delta p= -\frac32 \lambda^2 + \frac{\bar{\o} ^2}{2},
  \ee
  which is also obtained by taking subtraction (\ref{1.14})-(\ref{1.15}).  The vorticity equation is written as
\bb\label{1.17}
\frac{\bar{D} \bar{\o} }{Dt}=\lambda \bar{\o},
\ee
which can be solved as
$$
\bar{\o} (X_3 (a,t),t) =\bar{\o}_0 (a) \exp\left[\int_0 ^t \lambda (X_3(a,s),s)ds\right]
$$
along the trajectory. This implies that $\bar{\o}(X_3 (a,t),t)=0$ for $a\in \mathcal{S}$ as long as classical solution persists. Hence, for $a\in \mathcal{S}$ (\ref{1.15}) can be written as
\bq\label{1.17a}
\frac{\partial\lambda(X_3 (a,t),t)}{\partial t}&=&\frac{\lambda^2(X_3 (a,t),t) }{2}+2\partial_{r}^2 p (X_3 (a,t),t)\n\\
&\geq& \frac{\lambda^2(X_3 (a,t),t) }{2}\qquad \forall t\in (0,T_1(a)).
\eq
The differential inequality (\ref{1.17a}) can be solved immediately to yield
$$
\lambda(X_3 (a,t),t)\geq \frac{2\lambda_0 (a)}{2-t\lambda_0 (a)}\quad
\forall t\in (0, T_*)\quad \mbox{with}\quad T_*=T_*(a):=\min\{ T_1 (a), \frac{2}{\lambda_0 (a)} \},
$$
which shows that $T_1(a) \geq 2/\lambda_0 (a)$ is not consistent with the fact that classical solution
persists until $T(a)$.
$\square$\\
\ \\

\noindent{\bf Proof of Theorem 1.2 }
By the hypothesis  the equation (\ref{1.14}) together with (\ref{1.17}) reduces to
\bb\label{1.19b}
\left\{\aligned
&\frac{\bar{D}\lambda}{Dt}=\frac{\lambda^2}{2}-\frac{\bar{\o} ^2}{2},\\
&\frac{\bar{D}\bar{\o}}{Dt}=\lambda \bar{\o}.
\endaligned \right.
\ee
This is exactly the same system studied by Constantin-Lax-Majda in \cite{con} with the material derivative replacing the partial derivative in time, which was proposed as a one dimensional model equation for the 3D Euler equations in the vorticity formulation.
Similarly to \cite{con} we set $\Theta=\lambda+i \bar{\o}.$ Then (\ref{1.19b})  becomes the following complex Riccati equation along the trajectory,
$$
\frac{\bar{D}\Theta}{Dt}=\frac{\Theta^2}{2},
$$
which can be solved explicitly as
\bb\label{1.20}
\Theta (X_3(a,t),t)=\frac{2\Theta_0 (a)}{2-\Theta_0 (a)t}=\frac{2\lambda_0 (a) +2i \bar{\o}_0 (a)}{2-[\lambda_0 (a)+i \bar{\o}_0 (a)]t}.
\ee
Taking imaginary and real parts of (\ref{1.20}) we obtain (\ref{1.19})-(\ref{1.19a}). $\square$\\
\ \\
\noindent{\bf Remark after the proof } In \cite{hou} Hou-Li also obtained a system of equations similar in form to
(\ref{1.19b}), but for a different pair of unknown functions under completely different assumptions. In our case the system is derived rigorously from the axisymmetric 3D Euler equation by taking the limit $r\to 0$, and assuming only
$\partial^2_r p=0$ on the $x_3-$axis.
\[ \mbox{\bf Acknowledgements} \]
The author would like to thank to  P. Constantin, Y. Brenier and F. Gallaire for helpful discussions and comments.  This research was supported partially by KRF Grant(MOEHRD, Basic Research Promotion Fund).

  \end{document}